\documentclass[11pt]{amsart}

\usepackage[psamsfonts]{eucal}
\usepackage{graphicx}
\usepackage{amsmath, amssymb, amsfonts}
\usepackage{fullpage}
\usepackage{setspace}
\usepackage{enumerate}

\newtheorem{theorem}{Theorem}[section]
\newtheorem{proposition}[theorem]{Proposition}
\newtheorem{lemma}[theorem]{Lemma}

\newtheorem{remark}[theorem]{Remark}

\def\R{{\mathbb R}}
\def\Q{{\mathbb Q}}

\def\Z{{\mathbb Z}}

\def\N{{\mathbb N}}

\def\ve{{\varepsilon}}

\newcommand{\dd}{\,\mathrm{d}}
\newcommand{\D}{\displaystyle}
\newcommand{\tb}[1]{\textbf{#1}}

\begin{document}
\title{Orthogonal Exponentials for Bernoulli Iterated Function Systems}

%\titlerunning{Orthogonal Exponentials for Bernoulli IFSs}
%\date{\today}

% Use \titlerunning{Short Title} for an abbreviated version of
% your contribution title if the original one is too long
\author[P.E.T. Jorgensen, K. Kornelson, K. Shuman]{ Palle E. T. Jorgensen, Keri Kornelson, and Karen Shuman}

\address{(\textrm{P.E.T. Jorgensen}) Department of Mathematics, University of Iowa, Iowa City, IA
52242} \email{jorgen@math.uiowa.edu}

\address{(\textrm{K. Kornelson}) Department of Mathematics \& Statistics, Grinnell College,
Grinnell IA 50112-1690} \email{kornelso@math.grinnell.edu}

\address{(\textrm{K. Shuman}) Department of Mathematics \& Statistics, Grinnell College,
Grinnell IA 50112-1690} \email{shumank@math.grinnell.edu}

%\thanks{ }%

\begin{abstract} We investigate certain spectral properties of the Bernoulli convolution measures on attractor sets arising
from iterated function systems on $\R$.  In particular, we examine collections of orthogonal
exponential functions in the Hilbert space of square integrable functions on the attractor.  We
carefully examine a test case $\lambda = \frac34$ in which the IFS has overlap.  We also determine
rational $\lambda = \frac{a}{b}$ for which infinite sets of orthogonal exponentials exist.
\end{abstract}

\maketitle

\begin{center}
\small \textit{Dedicated with respect and fondness to Larry Baggett.} \normalsize \end{center}

\renewcommand{\baselinestretch}{1.2}\large \normalsize

\section{Introduction}

This work examines the spectral properties of a class of measures on fractals which arise from
affine iterated function systems (IFSs) on $\R$. The fractals $X$ are compact subsets of $\R$,
which may or may not have Hausdorff dimension less than $1$.  Even though such sets do not have a
group structure or Haar measure, we are able to identify a substitute.  Associated with each
compact set $X$ is a measure $\nu$, often called a Hutchinson measure, with reference to
\cite{Hut81}.  The measure $\nu$ is a probability measure having support $X$.  By analogy to a Haar
measure, $\nu$ is uniquely determined by the maps which characterize the iterated function system
and exhibits an invariance property arising from these maps.

The Bernoulli affine IFS on $\R$ is given by two functions
$$\tau_{b_1}(x)= \lambda(x-b_1), \qquad \tau_{b_2}(x) = \lambda(x+b_2),$$ where $B = \{b_1, b_2\}$
and the parameter $\lambda \in (0,1)$.  The measures $\nu$ arising from a Bernoulli IFS are called
Bernoulli convolution measures. These measures have a long history which can be studied, for example, in \cite{Ka85}.

We will examine orthogonality of exponential functions in the Hilbert space $L^2(X,\nu)$, in order
to better understand the Bernoulli convolution measures.  Jorgensen and Pedersen
\cite{JP1996, JP1998} recently found examples of Bernoulli measures for which the Hilbert space
$L^2(X,\nu)$ has an orthonormal basis of exponential functions.  In such cases, we say that $\nu$
is a spectral measure.

These known examples have certain properties in common.  The Bernoulli IFS maps are of the form
$\tau_{b_1}(x)= \lambda(x-b_1)$ and $ \tau_{b_2}(x) = \lambda(x+b_2),$ where $b_1$ and $ b_2$ are
integers, and $\lambda = \frac{1}{n}$ is the reciprocal of a natural number.  These conditions on
$\lambda$ and $B$ describe what we will call the \textit{rational case} of a Bernoulli affine IFS.
In recent papers \cite{DJ2006b,JP1998}, an additional condition to the rational case -- Hadamard
duality -- is discovered which will guarantee that an orthonormal Fourier basis exists
corresponding to the rational cases of $\lambda$ and $B$.

In this paper, we consider the non-rational cases for which the scaling factor $\lambda$ is not
thus restricted.  In such examples, we also find that the values $b_1$ and $b_2$ are not integers.
The arguments from dynamics and random-walk theory which can be used in the rational cases no
longer apply, so the problem of determining whether Fourier bases exist, and constructing them if
they do, are much harder.

\subsection{Motivation}
A leading theme in harmonic analysis in general \cite{Fol95}, and
in the work of Larry Baggett and co-authors in particular is that of building basis decompositions,
or direct integral decompositions, out of geometric composite structures.  Such processes generally
involve a group structure, as in the case of semidirect product groups, and the use of induced
representation constructions, both for discrete groups and continuous groups. This setting
accommodates a variety of applications, including wavelet analysis which relates to $ax + b$- like
groups and time-frequency analysis which involves Heisenberg-like groups.  For examples of results
in these areas, see \cite{BCMA91, BCMO95, BMM99}.

        There is recent evidence \cite{Jor06a, DuJo06b} that the question of basis decompositions can be examined on a class
of fractals, which are attractors of so called affine Iterated Function Systems (IFSs), including
Cantor's middle-third example; see also \cite{Hut81}. As it turns out, a number of the classical
tools from the group case still work, but there is no longer a group structure to work with.

Without the groups, we must first look for a substitute for Haar measure. Thanks to a construction
of Hutchinson \cite{Hut81}, this is available in the form of equilibrium measures $\nu$ with
support on the fractal attractor set $X$ of the IFS.  We can therefore describe a precise notion of
an orthogonal Fourier expansion in the Hilbert space $L^2(X,\nu)$.  It turns out in fact that the
equilibrium measures will typically be singular with respect to Lebesgue measure; see
\cite{JKS07}.

          While there is already some work on expansion problems in the context of affine IFS fractals \cite{DuJo06b,
JP1996}, in these papers the class of admissible IFSs with equilibrium measures $\nu$ is restricted
in several ways such that the Hilbert space $L^2(X,\nu)$ will carry an orthogonal Fourier expansion
based on complex exponentials.  However the systems considered so far have a rather restricted and
rigid Diophantine structure, and they do not admit overlap. They also do not admit continuous
deformations.  The affine IFSs considered here are not restricted to the conditions imposed in the
earlier studies \cite{DuJo06b, JP1996}.

\subsection{Related Topics}
 The analysis of fractal measures is motivated by and has influences on a variety of areas
outside fractal theory itself.  Some such fields include tiling spaces (see e.g., \cite{BaDi05}),
lacunary expansions (\cite{Tri04, Ko03}), and random Fourier series (\cite{Ka85, MaPi81}) . Each of
these subfields offers a variation of the general theme of recursive constructions based on some
notion of self similarity.  Moreover, each area invites an approach which mixes tools from operator
algebra theory and from probability. Together this part of mathematics stands at the crossroads of
operator algebras, basis constructions, and dynamics.  This study also interacts with parallel
developments in the study of wavelets (including wavelets on fractals \cite{DuJo06a, DuJo06b,
Jor06a}), and of iterated function systems; see e.g., \cite{DuJo06c}.

We summarize here a few relevant facts to illustrate the connections and common themes to these diverse topics.

\begin{enumerate}

\item Tiling spaces \cite{BaDi05} are important for our
understanding of diffraction in molecular structures that form quasi crystals. The simplest tiles
in $\R^d$ come from translations by rank $d$ lattices; and others involve both translations and
matrix operations. A pioneering paper by Fuglede \cite{Fu74} suggested a close connection between
the Fourier bases and translation tiles, and this connection was clarified in later works such as
\cite{Ta04}.

\item Lacunary Fourier series on the line refers to Fourier
expansions where there are infinitely many gaps in the frequency variable which separate powers of
a selected finite set of numbers. By inspection, one checks that the ``fractals in the large''
described in \cite{JP1996} which arise in the known examples as the spectra of affine IFSs fit this
pattern.

\item Random series of functions (as in  \cite{Ka85,MaPi81})
naturally generalize the first affine IFS systems -- those constructed as infinite
convolutions of independent Bernoulli variables as in \cite{E1939}. The notion
``random'' here refers to the study of Fourier series where the Fourier coefficients are
random variables.

\end{enumerate}

\subsection{Overview}
In Section \ref{sec:defs}, we carefully describe the Bernoulli affine iterated function system for
parameter $\lambda \in (0,1)$ and the resulting attractor set $X_{\lambda}$ and Hutchinson measure
$\nu_{\lambda}$.  We also describe some of the early results about this measure, by Erd\H{o}s and
others.  Next, we establish the notation and elementary results which will be used in the following
sections.  In Section \ref{sec:3/4} we determine that there do exist infinite collections of
orthogonal exponentials for the special case where $\lambda = \frac34$.  We find that some of these
collections are maximally orthogonal, in the sense that such a collection is not properly contained
in another orthogonal collection.

In Section \ref{sec:rational}, we determine for which rational values of $\lambda$ there exist
infinite orthogonal collections of orthogonal exponentials. Our main result is the following, which
is stated in the paper as Theorem \ref{thm:generallambda}. \textit{Given $\lambda = \frac{a}{b}$,
if $b$ is even, then there exist infinite families of orthogonal exponentials in the corresponding
Hilbert space.  If $b$ is odd, then every collection of mutually orthogonal exponentials must be
finite.}

In Section \ref{sec:total}, we conjecture that none of these collections of orthogonal exponentials
for $\lambda = \frac34$ are actually orthonormal bases for the Hilbert space.  This conjecture is
based on numerical evidence.  It would support a conjecture in \cite{DJ2006b} that no orthonormal
bases exist in the non-rational cases.

%%%%%%%%%%%%%%%%%%%%%%%%%%%%%%%%%%%%%%%%%%%%%%%%%%%%%%%%%%%%%
\section{Affine IFSs and their associated invariant measures}\label{sec:defs}
%%%%%%%%%%%%%%%%%%%%%%%%%%%%%%%%%%%%%%%%%%%%%%%%%%%%%%%%%%%%%%

We study families of exponentials which are mutually orthogonal with respect to invariant measures
$\nu_{\lambda}$ associated with affine iterated function system (IFS) on the real line.  As in
\cite{DJ2006a}, we will consider affine IFSs which are determined by two affine maps on $\R$
\begin{equation}\label{tauplusminus}
\tau_{+}(x):=\lambda x + 1
 \qquad\text{ and}\qquad
\tau_{-}(x):=\lambda x-1,
\end{equation}%Example 3.2 in DJ2006a
where $\lambda \in (0,1)$.

Because $\lambda = \frac{3}{4}$ arises as the first value in a special family of values of
$\lambda$ considered by Dutkay and Jorgensen (\cite[Theorem 4.5]{DJ2006a}), we first explore
families of exponentials associated with the affine IFS with parameter $\lambda = \frac{3}{4}$.  It turns out
that this value of $\lambda$ provides good intuition for what happens in the more general case
$\lambda \in (\frac{1}{2}, 1)\cap\Q$.  We then generalize our results for $\lambda = \frac{3}{4}$
to $\lambda = \frac{a}{b}$, where $a,b\in\N$.

When two affine maps as in (\ref{tauplusminus}) generate an IFS, the resulting measure is called an
\textit{infinite Bernoulli convolution measure} (Lemma \ref{lem:convmeas}). Infinite Bernoulli
convolution measures are characterized by the infinite product structure of their Fourier
transforms.  To be precise, when $\lambda\in (0,1)$, the product $\prod_{n = 0}^{\infty} \cos(2\pi
\lambda^n t)$ is the Fourier transform of an infinite Bernoulli convolution measure.  

The term ``infinite convolution'' is not mysterious at all, since the Fourier transform converts
convolutions to products.  However, one might ask why these measures have the name ``Bernoulli''
attached to them.   These measures arise in the work of Erd\H{o}s and others via the study of the
random geometric series $\sum \pm \lambda^n$ for $\lambda \in (0,1)$, where the signs are the
outcome of a sequence of independent Bernoulli trials.  In other words, we could consider the signs
to be determined by a string of fair coin tosses.  This makes $\sum \pm \lambda^n$ a random
variable, i.e. a measurable function from a probability space into the real numbers $\R$.  In Erd\H{o}s's language, $\beta$ is the distribution of the random variable $X$ is defined on the probability space $\Omega$ of all infinite sequences of $\pm 1$.  The measure on $\Omega$ is the infinite-product measure resulting from assigning $\pm1$ equal probability $\frac{1}{2}$.  The random variable $X$ takes on a specific real value for each sequence from $\Omega$.  This distribution $\beta$, then, is the
familiar Bernoulli distribution from elementary probability theory. The infinite Bernoulli
convolution measure is determined by the distribution $D$ of the random variable $\sum \pm
\lambda^n$, which can be constructed from infinite convolution of dilates of $\beta$:
\begin{equation*}
D := \beta(x) * \beta(\lambda^{-1}x) * \beta(\lambda^{-2}x)*\ldots*\beta(\lambda^{-n}x)*\ldots.
\end{equation*}

These Bernoulli convolution measures have been studied from at least the mid-1930s in various
contexts.  There seems to have been a flurry of activity in the 1930s and 1940s surrounding these
measures.  Jessen and Wintner study these measures in their study of the Riemann zeta function in
their 1935 paper \cite{JW1935}; in 1939 and 1940, Erd\H{o}s published two important papers about
these measures \cite{E1939, E1940}.  In the 1939 paper, Erd\H{o}s proved that if $\alpha$ is a
Pisot number (that is, $\alpha$ is a real algebraic integer greater than $1$ all of whose
conjugates $\tilde{\alpha}$ satisfy $|\tilde{\alpha}| < 1$), then the infinite Bernoulli
convolution measure associated with $\lambda = \alpha^{-1}$ is singular with respect to Lebesgue
measure.  However, more recently, Solomyak proved that for almost every $\lambda \in (\frac{1}{2},
1)$, the measure $\nu_{\lambda}$ is absolutely continuous with respect to Lebesgue measure
\cite{Sol95}.

Bernoulli convolution measures arise in the study of affine iterated function systems because they
possess a special \textit{invariance} property---a property which is described in Hutchinson's 1981
theorem \cite[Theorem 2, p. 714]{Hut81}.  In the cases studied here, the invariance property of the
measure can be written
\begin{equation}\label{eqn:invariance}
\nu_{\lambda} = \frac{1}{2}(\nu_\lambda\circ \tau_{+}^{-1} + \nu_\lambda\circ
\tau_{-}^{-1}).
\end{equation}
We will now see that the measure $\nu_\lambda$ satisfying this invariance property is a Bernoulli
convolution measure.  We note that the following lemma is not a new result; we state and sketch a
proof of the lemma to show the connection between affine IFSs and Bernoulli convolution measures.
See \cite[Lemma 4.7]{DJ2006a}.

\begin{lemma}\label{lem:convmeas}Given a fixed $\lambda \in (0,1)$, the Fourier transform of $\nu_{\lambda}$ satisfies the equation
\begin{equation}\label{eqn:product}
\widehat{\nu_{\lambda}}(t) =
 \prod_{n = 0}^{\infty}
\cos\Bigl(2 \pi \lambda^n t \Bigr).
\end{equation}
That is, the measure arising from the IFS (\ref{tauplusminus}) is a Bernoulli convolution measure.
\end{lemma}

\begin{proof}We use Hutchinson's invariance property in Equation (\ref{eqn:invariance}).  This invariance property has integral form
\begin{equation}\label{eqn:intinvariance}  \int f(x) \dd \nu_{\lambda}(x) = \frac{1}{2}
\left[ \int f(\tau_-(x)) \dd \nu_{\lambda}(x) + \int f(\tau_+(x)) \dd \nu_{\lambda}(x)\right].
\end{equation}  We use this relation to show that the Fourier transform of $\nu_{\lambda}$ has the structure of an infinite Bernoulli convolution measure.
\begin{eqnarray*}
 \widehat{\nu_{\lambda}}(t) &=& \int e^{2\pi i x t} \dd \nu_{\lambda}(x) \\&=& \frac12 \int e^{2\pi
i (\lambda x-1) t} + e^{2\pi i (\lambda x+1) t} \dd \nu_{\lambda}(x) \quad \mbox{by Eq. (\ref{eqn:intinvariance})} \\
&=& \frac12 (e^{-2\pi i  t}+e^{2\pi i t})\int e^{2 \pi i \lambda x t} \dd \mu_{\lambda}(x)\\
&=&  \cos(2\pi t) \widehat{\nu_{\lambda}}(\lambda t)\\ &=&  \cos( 2\pi t)\cos(2\pi \lambda t)
\widehat{\nu_{\lambda}}(\lambda^2 t) \\ & \vdots & \qquad \vdots
\\ &=&  \left( \lim_{n \rightarrow \infty} \widehat{\nu_{\lambda}}(\lambda^nt)\right)\prod_{n=0}^{\infty} \cos(2\pi \lambda^n t)
 \end{eqnarray*}

The limit $\lim_{n \rightarrow \infty} \widehat{\nu_{\lambda}}(\lambda^nt)$ exists and is equal to
$\widehat{\nu_{\lambda}}(0)=1$, since it is known that $\nu_{\lambda}$ has no atoms and is a
probability measure.  Dutkay and Jorgensen show that $\nu_{\lambda}$ has no atoms in \cite[Corollary 6.6 ]{DJ2006a}; the fact that $\nu_{\lambda}$ is a probability measure follows from Hutchinson's theorem \cite[Section 4.4]{Hut81}. \qed
\end{proof}

Let $X_{\lambda}$ be the support of $\nu_{\lambda}$---that is, $X_{\lambda}\subset \R$ is the attractor for the IFS (\ref{tauplusminus}).
 For $\Lambda \subset \R$, we determine conditions under which the set of exponentials $\{e^{2\pi i \ell
\cdot}\,:\, \ell \in \Lambda\}$ is an orthonormal basis in the Hilbert space $L^2(X_{\lambda},
\nu_{\lambda})$.  In other words, we explore whether $\Lambda$ is a spectral set for the Hilbert
space.

It was shown in \cite{JP1998} that when $\lambda = \frac13$, there is no orthonormal basis of
exponentials, but when $\lambda = \frac14$ there is such an ONB.  In this paper, we are
particularly interested in values of $\lambda >\frac12 $, for which the IFS has overlap.  After some
general formulas, we will examine the special case $\lambda = \frac34$.   In Section
\ref{sec:rational}, we state some generalizations to other rational values of $\lambda$.

We wish to show that the infinite product in Equation (\ref{eqn:product}) is zero if and only if
one of the factors is zero. In the following lemma, we have omitted the $2\pi$ for notational
convenience.

\begin{lemma}\label{product_lemma}  Suppose $t$ is a fixed real number and that $\lambda \in(0,1)$.  There exists $N\in\N$ and $c > 0$ such that
\begin{equation}
\prod_{n = N}^{\infty} \cos(\lambda^n t)\geq c.
\end{equation}
In other words, if for some $t_0\in\R$,
\[ \prod_{n = 0}^{\infty}\cos(\lambda^n t_0)= 0,\]
then one of the factors of the product must be $0$.
\end{lemma}

\tb{Proof:  }To start, we note that
\begin{equation*}
\prod_{n = N}^{\infty} \cos(\lambda^n t)\geq c
\Leftrightarrow
\sum_{n=N}^{\infty} \ln\Bigl(\cos(\lambda^n t)\Bigr)\geq \ln(c).
\end{equation*}
The Taylor expansion of cosine around $0$ yields
\begin{equation*}
 \cos(\lambda^n t)
= 1 - \sum_{k = 1}^{\infty} (-1)^{k+1}\frac{\lambda^{2kn}}{(2k)!}t^{2k}.
\end{equation*}
Define
\[ \ve_{n}(t):= \sum_{k = 1}^{\infty} (-1)^{k+1}\frac{\lambda^{2kn}}{ (2k)!}t^{2k};\]
since $-1 \leq 1 -\ve_{n}(t)\leq 1$ for all $n\in\N$,  we know that $\ve_n(t)\geq 0$ for all $n\in\N$.

We can choose $N_1$ such that for all $n > N_1$, \; $\displaystyle \ve_n(t) \leq \frac{\lambda^{2n}}{2}t^2$, and then we can choose $N_2$ such that for all $n>N_2$,\; $\displaystyle\frac{\lambda^{2n}}{2}t^2<1$.  We now consider the Taylor expansion of $\ln(1-\ve_n(t))$, which is valid when $|\ve_n(t)|<1$:
\begin{equation*}
\ln\Bigl(\cos(\lambda^n t)\Bigr)
= \ln(1 - \ve_n(t))
= -\sum_{k =1}^{\infty}\frac{\ve_n(t)^k}{k} .
\end{equation*}
Finally, we choose $N_3$ such that for all $n > N_3$,\;
$\displaystyle \sum_{k = 1}^{\infty}\frac{\ve_n(t)^k}{k} \leq 2 \ve_n(t)$.

Now, for $N > \max\{N_1, N_2, N_3\}$, we have
\begin{eqnarray*}
\sum_{n=N}^{\infty} \ln\Bigl(\cos(\lambda^n t)\Bigr) &=& \sum_{n=N}^{\infty}\Biggl(-\sum_{k = 1}^{\infty}\frac{\ve_n(t)^k}{k}\Biggr)\\&\geq& \sum_{n=N}^{\infty} -2 \ve_n(t) \\&\geq& \sum_{n=N}^{\infty} -\lambda^{2n}t^2.
\end{eqnarray*}
But $\sum_{n=N}^{\infty} \lambda^{2n}$ is a convergent geometric series, so set
\begin{equation*}
c = \exp\Bigl(-t^2\sum_{n=N}^{\infty}\lambda^{2n}\Bigr).
\end{equation*}
We have now found $N$ and $c$ such that
\begin{equation*}
\prod_{n = N}^{\infty} \cos(\lambda^n t)\geq c.
\end{equation*}
\qed

Our goal in this paper is to study collections of orthogonal exponentials in the Hilbert space $L^2(X_{\lambda}, \nu_{\lambda})$.   We now observe that the Fourier transform $\widehat{\nu_{\lambda}}$ of the Bernoulli measure
 arises in the inner product of exponential functions in the Hilbert space
$L^2(X_{\lambda}, \nu_{\lambda})$.
\begin{equation}\label{orthognality_condition}
\langle e^{2\pi i \ell t}, e^{2\pi i \ell' t} \rangle = \int e^{2\pi i \ell t}\overline{e^{2\pi i
\ell' t}}\: d\nu_{\lambda}(t) = \widehat{\nu_{\lambda}}(\ell-\ell')
\end{equation}

Using Lemmas \ref{lem:convmeas} and  \ref{product_lemma}, we can conclude that the exponential functions $e_{\ell}$ and $e_{\ell'}$ are orthogonal if and only if $\cos(2\pi \lambda^n(\ell-\ell')) = 0$ for some $n \in \N_0$.  (To simplify notation, we will use $e_{\ell}$ to denote the exponential function $e^{2\pi i \ell \cdot}$.)
%%%%%%%%%%%%%%%%%%%%%%%%%%%%%%%%%%%%%%%%%%%%%%%%%%%%%%%%%%%%%%%%%%%%%
\section{Orthogonal exponentials with respect to $\nu_{\frac{3}{4}}$}\label{sec:3/4}
%%%%%%%%%%%%%%%%%%%%%%%%%%%%%%%%%%%%%%%%%%%%%%%%%%%%%%%%%%%%%%%%

In order to understand the issues involved with the study of orthogonal exponentials, and the possible existence of orthonormal bases, we study the
special case $\lambda = \frac34$.

\begin{lemma}\label{lem:orthog}  The function $\widehat{\nu_{\frac34}}(t)$ is equal to zero if and only if $t\in \frac{4^{n-1}}{3^n}(1 + 2\Z)$ for some $n\in\N_0$.  \end{lemma}

\begin{proof}
Choose $n\in\N_0$ and $\ell\in\Z$, and let $t = \frac{4^{n-1}}{3^n}(2\ell+1)$.  Then
\[\cos\Bigl(2 \pi \Bigl(\frac{3}{4}\Bigr)^n t \Bigr)=  \cos\Bigl(2 \pi \frac{3^n}{4^n} \cdot \frac{4^{n-1}}{3^n}(2\ell+1)\Bigr) =
\cos\Bigl(\frac{\pi}{2}(2\ell+1)\Bigr) = 0.\]  Since this is one of the factors in the infinite product \[ \widehat{\nu_{\frac34}}(t)=  \prod_{n = 0}^{\infty}
\cos\Bigl(2 \pi \Bigl(\frac34\Bigr)^n t \Bigr),\]
we see that  $\widehat{\nu_{\frac34}}(t)$ is $0$.

Conversely, suppose that for some $t_0 \in \R$ that $\widehat{\nu_{\frac34}}(t_0)=0$.  By Lemma \ref{product_lemma}, this implies that for some  $n\in \N_0$, that the factor  $\cos\Bigl(2 \pi \Bigl(\frac{3}{4}\Bigr)^n t_0\Bigr)$ is equal to zero. 
 Therefore, the quantity $2 \pi \Bigl(\frac{3}{4}\Bigr)^n t_0$ must be an odd multiple of $\frac{\pi}{2}$, so there 
 exists $\ell\in\Z$ such that
\[2 \pi \Bigl(\frac{3}{4}\Bigr)^n t_0 = \frac{\pi}{2}(2\ell +1),\]
which then gives
\[ t_0 = \frac{4^{n-1}}{3^n}(2\ell+1).\]

As a result, we have established that $\widehat{\nu_{\frac34}}(t_0) = 0$ if and only if one of the
terms in the infinite product defining $\widehat{\nu_{\frac34}}$ is $0$.  In other words, \[
\widehat{\nu_{\frac34}}(t_0) = 0 \Leftrightarrow t_0 \in \Biggl\{\frac{4^{n-1}}{3^n}(1 + 2\Z) :
n\in\N\Biggr\}.\] \qed
\end{proof}

We will begin our examination of the sets of orthogonal exponentials for $L^2(X_{\frac34},\nu_{\frac34})$ by finding sets $\Lambda \subset \R$ such that $\{ e_{\ell} : \ell\in\Lambda\}$ is a
mutually orthogonal set of functions with respect to $\nu_{\frac34}$.  By the preceding discussion, we have the following proposition.

\begin{proposition}\label{prop:orthogonalset}  The exponential functions $\{ e_{\ell} : \ell\in\Lambda\}$ are pairwise orthogonal if and only if for $\ell, \ell' \in \Lambda$, we have $\ell-\ell' \in \mathcal{O}$, where we define
\begin{equation}\label{eqn:orthogonaltest} \mathcal{O}= \Biggl\{\frac{4^{n-1}}{3^n}(1 + 2\Z) : n\in\N_0 \Biggr\}.
\end{equation}

\end{proposition}
 \begin{proof}  This follows from Equation (\ref{orthognality_condition}), Lemma \ref{product_lemma}, and Lemma \ref{lem:orthog}. \qed \end{proof}

\begin{theorem}\label{lambda_set_theorem}
There exist infinitely many infinite sets $\Lambda$ such that $\{ e_{\ell} :
\ell\in\Lambda\}$ is a mutually orthogonal set of functions with respect to the measure
$\nu_{\frac34}$.
\end{theorem}

\begin{proof} Define the set $\Lambda_k$ for each $k \in \N$ as follows.
\begin{equation}\label{lambda_k_defn}
\Lambda_k := \Biggl\{\frac{4^j}{3^k}: j\in\N, j \geq k-1 \Biggr\}\cup\{0\}
\end{equation}
If $\ell = \frac{4^p}{3^k}$ and $\ell' = \frac{4^q}{3^k}$ ($p,q\geq k-1$, $p > q$) are two non-zero
elements of $\Lambda_k$, then
\begin{equation}
\ell - \ell' = \frac{4^p}{3^k} - \frac{4^q}{3^k} = \frac{4^q}{3^k}(4^{p-q}-1),
\end{equation}
and $4^{p-q}-1$ is an odd integer.  Now, multiply through by $\frac{3^{q-k+1}}{3^{q-k+1}}$:
\begin{equation}
\ell - \ell' =  \frac{4^q}{3^k}\frac{3^{q-k+1}}{3^{q-k+1}} (4^{p-q}-1)= \frac{4^q}{3^{q+1}}(3^{q-k+1})(4^{p-q}-1),
\end{equation}
and $(3^{q-k+1})(4^{p-q}-1)$ is still an odd integer.  Therefore, $\widehat{\nu_{\frac34}}(\ell - \ell') = 0$.

If $\ell \neq 0$ and $\ell' = 0$, then
\[ \frac{4^p}{3^k} = \frac{4^p}{3^k}\cdot \frac{3^{p-k+1}}{3^{p-k+1}} = \frac{4^p}{3^{p+1}} 3^{p-k+1},\] \qed
and $3^{p-k+1}$ is an odd integer since $p \geq k-1$.
\end{proof}

If we draw a diagonal diagram of the $\Lambda_k$ sets in Theorem \ref{lambda_set_theorem}, we have

\[\begin{matrix}
\mathbf{j:  }                                    & & 0         & 1      &2         &3          &4            &\cdots \\
\mathbf{\Lambda_1:  } & 0 & \frac{1}{3} & \frac{4}{3} & \frac{16}{3} & \frac{64}{3} & \frac{128}{3} & \cdots\\
\mathbf{\Lambda_2:  } & 0 &             & \frac{4}{9} & \frac{16}{9} & \frac{64}{9} & \frac{128}{9} & \cdots\\
\mathbf{\Lambda_3:  } & 0 &             &             & \frac{16}{27} & \frac{64}{27} & \frac{128}{27} & \cdots\\
\vdots                & 0 &             &             &               & \ddots        &                &\\
\end{matrix} \]

These sets are certainly not the only infinite sets satisfying the condition of Proposition \ref{prop:orthogonalset}.  Since our condition tests the differences between elements, any of the above can be translated by a real number $\alpha$.  We choose representative sets by making the requirement that $0$ be an element of each $\Lambda_k$.   Similarly, if every element is a $\Lambda_k$ set is multiplied by the same odd integer, the differences remain elements of the set $\mathcal{O}$.

One natural question is whether some of these $\Lambda_k$ sets can be combined to form larger
collections of orthogonal exponentials.  We find, however, that the orthogonality condition from
Proposition \ref{prop:orthogonalset} is lost if we take the union of different $\Lambda_k$ sets.

\begin{proposition}\label{adding_elements_proposition}
Suppose $\ell\neq 0$, $\ell \in\Lambda_k$ where $k > 1$.  Then the set $\Lambda_1\cup\{\ell\}$ does
not form a mutually orthogonal family of exponential functions.
\end{proposition}

\begin{proof}
Let $\ell = \frac{4^j}{3^k}$, where $k>1$ and $j \geq k$.  We can show that $e_{\ell}$ and
$e_{1/3}$ are not orthogonal, and therefore we cannot add $\ell$ to the set $\Lambda_1$ to build a
larger family of orthogonal exponential functions.

\begin{eqnarray*} \ell - \frac13 &=& \frac{4^j}{3^k}-\frac13 \\ &=& \frac{4^j-3^{k-1}}{3^k}
 \end{eqnarray*} The numerator in the last expression is still an integer since $k-1 >
0$, and we also observe that $4$ does not divide the numerator.  This means that if we are to write
this difference in the form of set $\mathcal{O}$ from (\ref{eqn:orthogonaltest}), our power of $4$
must be zero, and therefore the power of $3$ must be 1.

$$ \ell - \frac13 = \frac{4^j-3^{k-1}}{3^k} = \frac{4^0}{3^1}\,\left(\frac{4^j-3^{k-1}}{3^{k-1}}\right)$$

Since $j \geq k > 1$, the second fraction cannot be an integer, and therefore this difference is
not an element of the set $\mathcal{O}$.  This proves that the two exponentials $e_{\ell}$ and  $e_{1/3}$
 are not orthogonal. \qed
\end{proof}

A similar argument shows that nonzero elements from $\Lambda_k$ cannot be combined with $\Lambda_j$
for $j \neq k$ while maintaining orthogonality.

We can merge a finite number of $\Lambda_k$'s and still form an orthogonal set, as long as we are
willing to throw out finitely many terms.  For example, $(\Lambda_2\cup\Lambda_1)\setminus \{\frac{1}{3}\}$
forms an orthogonal set.  However, experimental evidence indicates that
$(\Lambda_2\cup\Lambda_1)\setminus \{\frac{1}{3}\}$ is not total (see Section \ref{sec:total}).

Rather than merging the $\Lambda_k$ sets, we can expand each of them to a larger collection of
orthogonal exponentials.

\begin{theorem}
Define the set $\Gamma_k$ for each $k \in \N$ as follows:
$$\Gamma_k = \left\{ \sum_{j=k-1}^p \frac{a_j 4^ j}{3^k}\,:\, p \,\mbox{ \em{finite}}, a_j \in \{0,1\}\right\} \bigcup \{0\}$$
Each set $\{e_{\gamma}\,:\, \gamma \in \Gamma_k\}$ is an orthonormal family in
$L^2(X_{\frac34},\nu_{\frac34})$.
\end{theorem}

\begin{proof}
We will demonstrate the proof for the set $\Gamma_1$.  The argument is similar for each $\Gamma_k$.

First, we show that each nonzero element in $\Gamma_1$ does belong to the set $\mathcal{O}$, which
proves that the differences with $0$ are in $\mathcal{O}$.  Let $a \in \Gamma_1, a \neq 0$.  Then
$\D a = \frac13 \sum_{i=0}^p a_i 4^i$, where $a_i \in \{0,1\}$.  Let $r$ be the smallest integer
such that $a_r \neq 0$, so $0 \leq r < p$.  We can then write $a$ in the form:

$$ a = \frac13 \sum_{i=r}^p a_i 4^i = \frac{4^r}{3} \sum_{i=r}^p a_i 4^{i-r} \left(\frac{3^r}{3^r}
\right) = \frac{4^r}{3^{r+1}} \left( 3^r \sum_{i=r}^p a_i 4^{i-r} \right)$$

Since $a_r =1$, we know that the sum in the last two expressions above is an odd integer, which
when multiplied by $3^r$ yields another odd integer.  Therefore, $a \in \mathcal{O}$.

Next, we must show that the difference of any two nonzero elements in $\Gamma_1$ is in
$\mathcal{O}$.  Let $a= \frac13 \sum_{i=0}^p a_i 4^i$ and $b = \frac13 \sum_{i=0}^q b_i 4^i$.  In
order to combine these, assume without loss of generality that $p \geq q$ and let $b_i=0$ for
$i=q+1, \ldots p$.  Therefore, we have $ a-b = \frac13 \sum_{i=0}^p (a_i-b_i)4^i $, where $a_i-b_i
\in \{-1,0,1\}$.  As above, let $r$ be the smallest integer such that $a_r - b_r \neq 0$, so $0
\leq r < p$. We can then write $a-b$ as follows:

$$ a-b = \frac13 \sum_{i=r}^p (a_i-b_i) 4^i = \frac{4^r}{3}\sum_{i=r}^p (a_i-b_i) 4^{i-r} \left(
\frac{3^r}{3^r} \right) = \frac{4^r}{3^{r+1}} \left( 3^r \sum_{i=r}^p (a_i-b_i) 4^{i-r} \right)$$

Since $a_r-b_r = \pm 1$, the sum in the last two expressions above is an odd integer, and so is its
product with $3^r$.  Therefore, $a-b \in \mathcal{O}$.

A parallel argument works for the other sets $\Gamma_k$, and in fact, demonstrates why each set
must start with powers of $4$ no more than one less than the power $k$ of $3$ in the denominator.
\qed
\end{proof}

For the same reasons that the sets $\Lambda_k$ cannot be combined while retaining orthogonality of
the exponentials (Proposition \ref{adding_elements_proposition}), the $\Gamma_k$ sets also cannot
be combined.   In fact, we find that the sets $\Gamma_k$ are maximal in a stronger sense as well.
For each $k$, $\Gamma_k$ is not strictly contained in another set for which all the exponentials
are pairwise orthogonal.  We will state the proof here for the set $\Gamma_1$ for ease of notation,
but remark that a parallel argument holds for each $\Gamma_k$.

\begin{theorem}\label{thm:maximal}  $\{e_{\gamma}\,:\, \gamma \in \Gamma_1\}$ is a maximally orthogonal
collection of exponentials for $L^2(X_{\frac34},\nu_{\frac34})$.  In other words, given $x \in \R
\setminus \Gamma_1$ there exists $\gamma \in \Gamma_1$ such that $e_x$ and $e_{\gamma}$ are not
orthogonal.
\end{theorem}

\begin{proof} First, note that since $0 \in \Gamma_1$, if $x \notin \mathcal{O}$, where $\mathcal{O}$
is the set given by Equation (\ref{eqn:orthogonaltest}) in Proposition \ref{prop:orthogonalset},
then $e_x$ is not orthogonal to $e_0$.  Thus, we can restrict to $x \in \mathcal{O}$, so $x =
\frac{4^{n-1}(2k+1)}{3^n}$ for some choice of $n \geq 1$ and $k \in \Z$.

\begin{itemize} \item[] \textit{Case $1$}, $(n=1).$\;  Take $x = \frac{2k+1}{3}$, but $x \notin \Gamma_1$.
Then $2k+1$ can be written in a base-$4$ expansion $2k+1 = \sum_{i=0}^p a_i 4^i$, where at least
one of $\{a_0, a_1, \ldots a_p\}$ is either $2$ or $3$.  Let $r$ be the smallest index such that
$a_r = 2$ or $3$.  If $a_r = 2$, let $\displaystyle \gamma = \frac{\sum_{i=0}^{r-1} a_i 4^i}{3}$ and if $a_r = 3$, let
$\displaystyle \gamma = \frac{(\sum_{i=0}^{r-1} a_i4^i)+ 4^r}{3}$.

In both cases above, we have $$ x-\gamma = \frac{2\cdot 4^r +\sum_{i=r+1}^p a_i 4^i}{3}.$$  Since
the numerator cannot be written as an odd multiple of a power of $4$, we have $x - \gamma \notin
\mathcal{O}$ which proves that the exponentials $e_x$ and $e_{\gamma}$ are not orthogonal by
Proposition \ref{prop:orthogonalset}.

\medskip

\item[]\textit{Case $2$}, $(n>1)$. \; Let $x = \frac{4^{n-1}(2k+1)}{3^n}$ where $n>1$.  It is possible that
powers of $3$ divide the odd integer $2k+1$, so we can cancel some of these if they exist to write
$x = \frac{4^{n-1}(2\ell+1)}{3^m}$ where either $m > 1$ and $2\ell+1$ is not divisible by $3$ or we
have $m=1$.

If $m >1$, then let $\gamma = \frac13$.  We find that $x-\gamma = \frac{4^{n-1}(2\ell+1) -
3^{m-1}}{3^m}$.  Since there is no way to cancel more powers of $3$, the numerator would need to be
divisible by $4^{m-1}$ if $x-\gamma$ were to be an element of $\mathcal{O}$.  We see, however, that
the numerator is not divisible by $4$.  Therefore the exponential functions $e_{\gamma}$ and
$e_{\frac13}$ are not orthogonal.

If $m =1$, then we have a situation similar to Case $1$.  Since $x = \frac{4^{n-1}(2\ell+1)}{3}$ but $x \notin \Gamma_1$ we have that
$2\ell+1$ has a base-$4$ expansion that includes at least one coefficient which is not a $0$ or a
$1$.  As above, given $2\ell+1 = \sum_{i=0}^p a_i4^i$, let $r$ be the smallest index for which $a_r
= 2$ or $3$.  If $a_r = 2$ let $\displaystyle \gamma = \frac{4^{n-1}\sum_{i=0}^{r-1}a_i4^i}{3}$ and if $a_r = 3$,
let $\displaystyle \gamma = \frac{4^{n-1}(\sum_{i=0}^{r-1}a_i4^i + 4^r)}{3}$.  This gives in both cases
$$x-\gamma = \frac{4^{n-1}(2\cdot4^r + a_{r+1}4^{r+1} + \cdots + a_p4^p)}{3}.$$  The numerator
cannot be expressed as an odd multiple of a power of $4$, so $x-\gamma \notin \mathcal{O}$ and
therefore, $e_\gamma$ and $e_x$ are not orthogonal.
\end{itemize}

This proves that the set $\Gamma_1$ is maximal, in the sense that there is no exponential function
$e_x, x \notin \Gamma_1$ which can be added to $\{e_{\gamma}\,:\, \gamma \in \Gamma_1\}$ to form an
orthogonal collection properly containing $\Gamma_1$. \qed
\end{proof}

%%%%%%%%%%%%%%%%%%%%%%%%%%%%%%%%%%%%%%%%%%%%%%%%%%%%%%%%%%%%%%%%%%%%
\section{Rational values of $\lambda$}\label{sec:rational}
%%%%%%%%%%%%%%%%%%%%%%%%%%%%%%%%%%%%%%%%%%%%%%%%%%%%%%%%%%%%%%%%%%%%
 We outlined in Section \ref{sec:defs} our rationale for focusing on the Bernoulli IFSs and on the
specific value $\lambda = \frac34$ for the scaling constant.  It is natural to next explore whether
these results are typical for other rational values of $\lambda$.  In this section, we find that
the orthogonality results from Section \ref{sec:3/4} do indeed extend to rational values of
$\lambda$ other than $\frac34$.

\begin{theorem}\label{thm:generallambda}  Let $\lambda \in \Q \cap (0,1)$ and let $\lambda = \frac{a}{b}$ be in reduced
form.  If $b$ is odd, then any collection of pairwise orthogonal exponential functions in the Hilbert space
$L^2(X_{\lambda}, \nu_{\lambda})$ can have only finitely many elements.  If $b$ is even, then there
exists a countably infinite collection of orthogonal exponentials in $L^2(X_{\lambda},
\nu_{\lambda})$.
\end{theorem}

Before we prove this theorem, we must find a new set $\mathcal{O}$ corresponding to the set
(\ref{eqn:orthogonaltest}) in Proposition  \ref{prop:orthogonalset} which will identify the zeros
of $\widehat{\nu_{\lambda}}(t)$ and thereby serve as a test for orthogonality of exponential
functions. We will assume that the collections under consideration all contain $0$.

\begin{lemma}  A set $\Gamma$  of real numbers containing $0$ has the property that  $\{e_{\gamma} \,:\, \gamma \in \Gamma \}$
is an orthogonal collection of exponentials in $L^2(X_{\lambda}, \nu_{\lambda})$ if and only if for
each $\gamma, \gamma' \in \Gamma$, we have $\gamma-\gamma' \in \mathcal{O}$, where we define

\begin{equation}\label{eqn:newOset}  \mathcal{O} := \left\{ \frac14
\Big(\frac{b}{a}\Big)^n (1+2\Z) \;:\; n \in \N_0 \right\}. \end{equation}
\end{lemma}

\begin{proof}
We showed in Section \ref{sec:defs} that $\langle e_{\gamma}, e_{\gamma'} \rangle =
\widehat{\nu_{\lambda}}(\gamma-\gamma') = \prod_{n=0}^{\infty} \cos(2\pi \lambda^n \gamma)$ for any
choice of $\lambda \in (0,1)$.  By Lemma \ref{product_lemma}, we know that
$\widehat{\nu_{\lambda}}(\gamma-\gamma')$ is zero if and only if for some $n$, the expression
$\cos(2 \pi \lambda^n (\gamma - \gamma'))$ is zero.  We reproduce the computations from Lemma
\ref{lem:orthog} to find our new set $\mathcal{O}$ corresponding to the set in Equation
\ref{eqn:orthogonaltest}.  If we take $\lambda= \frac{a}{b}$, we find that
$$2\pi \lambda^n (\gamma - \gamma') = (2k+1)\frac{\pi}{2}  \Leftrightarrow \gamma -\gamma' = \frac14
\frac{b^n}{a^n}(2k+1) \qquad k \in \Z.
$$

Since we are taking $0$ to always be an element in our sets of frequencies $\Gamma$, this means
each element of $\Gamma$ must already be an element of $\mathcal{O}$, and the differences between
any two elements must also be in $\mathcal{O}$.
\qed \end{proof}

\begin{proof}[Theorem \ref{thm:generallambda}]
Let $\lambda = \frac{a}{b}$, where $a$ and $b$ are relatively prime.  Let $\Gamma$ be a collection
of real numbers containing $0$ such that each nonzero element is in $\mathcal{O}$.  We can
accomplish much of this proof using parity arguments, so we consider individually the cases for
which $a,b$ are even/odd.

\begin{enumerate}
\item \textbf{$a$, $b$ odd:}\;    Let $\gamma, \gamma' \in \Gamma$ be
distinct and nonzero, so that $\gamma = \frac{b^n}{a^n}(2k+1)$ and $\gamma' =
\frac{b^{n'}}{a^{n'}}(2k'+1)$. Without loss of generality, let $n \geq n'$.

$$ \gamma - \gamma' = \frac14 \left[\frac{b^n(2k+1)-b^{n'}a^{n-n'}(2k'+1)}{a^n}\right] $$
Since the numerator above is an even number, it cannot be written as a product of $b^p$ and some
odd integer.  Multiplying this expression in both numerator and denominator by powers of $a$ won't
resolve this.  Therefore, we find that in the case where $a$ and $b$ are both odd, an orthogonal
set $\Gamma$ can only contain one nonzero frequency.

\item (\textbf{$a$ even, $b$ odd})\;  Let $\gamma \neq 0$ be an element of $\mathcal{O}$, so that
$\gamma = \frac14 (\frac{b}{a})^n (2k+1)$.  We will show that any collection $\Gamma$ containing
$0$ and $\gamma$ and having the property that the difference between any two elements is an element
of $\mathcal{O}$ must be a finite set.   Let $\gamma'\in \Gamma$, so that $\gamma'-0 = \frac14
(\frac{b}{a})^{n'}(2k'+1)$ in $\mathcal{O}$.

\begin{lemma}\label{lem:abcases}  Let $b$ be odd and $a$ be even.  Fix $\gamma \in \mathcal{O}$.
For every $\gamma' \in \mathcal{O}$ with $\gamma' \neq \gamma$, exactly one of the following is
true.
\begin{enumerate}\renewcommand{\labelenumi}{(\roman{enumi})}
\item $n=n'$ and $\gamma-\gamma' = \frac14 (\frac{b}{a})^p(2m+1)$ with $p<n$;
\item $n \neq n'$ and $\gamma-\gamma' = \frac14 (\frac{b}{a})^p(2m+1)$ with $p=\max\{n,n'\}$;
\item $\gamma - \gamma' \notin \mathcal{O}$.
\end{enumerate}

\end{lemma}

\begin{proof}
Let $n>n'$.  We have $$\gamma-\gamma' = \frac14 \left(\frac{b^n(2k+1)+b^{n'}a^{n-n'}(2k'+1)}{a^n}
\right). $$  Since $a$ is even, it does not divide the odd integer $b^n(2k+1)$, and therefore does
not divide the numerator.  We therefore have $\gamma - \gamma'$ in one of the following rational
forms:
$$ \frac14 \left(\frac{b}{a}\right)^p
\frac{b^{n-p}(2k+1)+b^{n'-p}a^{n-n'}(2k'+1)}{a^{n-p}} \qquad \mbox{for }\; p \leq n'$$ or
$$\frac14 \left(\frac{b}{a}\right)^p \frac{b^{n-n'}(2k+1)+a^{n-n'}(2k'+1)}{a^{n-p}b^{p-n'}}
\qquad \mbox{for }\;  n'< p<n $$ or $$\frac14 \left(\frac{b}{a}\right)^p
\frac{b^{n-n'}a^{p-n}(2k+1)+a^{p-n'}(2k'+1)}{b^{p-n'}} \qquad \mbox{for }\;  n'<n<p.$$

In the first two cases, the final factor in the expression cannot be an odd integer since $a$ does
not divide the numerator.  In the third case, the numerator of the final fraction is even while the
denominator is odd.  This fraction could be an integer but cannot be an odd integer.  Therefore,
the only remaining possibility is $p=n$.

An identical argument holds for $n<n'$.  Therefore, if $n \neq n'$, either $p=\max\{n,n'\}$ or
$\gamma-\gamma' \notin \mathcal{O}$.

Let $n=n'$.  Then $$\gamma-\gamma' =\frac14 \left(\frac{b}{a}\right)^n(2k-2k') $$ If $p \geq n$,
the rational form of $\gamma-\gamma'$ is $$ \frac14 \left(\frac{b}{a}\right)^p
\frac{2(k-k')a^{p-n}}{b^{p-n}}.$$  As above, since the numerator of the last factor is an even
integer and the denominator is odd, their ratio cannot be an odd integer.  Therefore, when $n=n'$,
either $p<n$ or $\gamma-\gamma' \notin \mathcal{O}$.  \qed
\end{proof}

Lemma \ref{lem:abcases} shows that given $\gamma \in \mathcal{O}$, we need only consider elements
$\gamma' \in \mathcal{O}$ which satisfy $n=n'$ or for which $n \neq n'$ and $p=\max\{n,n'\}$.  We
will show that we can only find finitely many such $\gamma'$ which also have their pairwise
differences in the set $\mathcal{O}$.

\begin{enumerate}
\item (\textit{$n=n'$})\; Suppose $\gamma' = \frac14 \left(\frac{b}{a}\right)^{n}(2k'+1)$.  We seek $\gamma'$
for which $$ \gamma-\gamma' = \frac14 \left(\frac{b}{a} \right)^{p}\frac{b^{n-p}2(k-k')}{a^{n-p}} ,
$$ where the final factor is an odd integer, and we have shown in Lemma \ref{lem:abcases} that
$p<n$ is a necessary condition.  Since $a$ and $b$ are relatively prime, we see that the expression
$\frac{b^{n-p}2(k-k')}{a^{n-p}}$ is an odd integer only if $\frac{2(k-k')}{a^{n-p}}$ is an odd
integer.

Suppose that $k$ and $k'$ are congruent modulo $a^{n+1}$.  Then $a^{n+1}$ divides $k-k'$, and
therefore $\frac{2(k-k')}{a^{n-p}}$ is an even integer.  Therefore, given our fixed $\gamma \in
\Gamma$, there can be only one $k'$ from each congruence class modulo $a^{n+1}$ to make up a
$\gamma'$ in $\Gamma$.

\item (\textit{$n > n'\geq 0$})\;   Since $p=n$ is our only option, by Lemma \ref{lem:abcases},
we seek $\gamma'$ such that
$$ \gamma-\gamma' = \frac14 \left(\frac{b}{a}\right)^n \left[(2k+1)-\frac{(2k'+1)a^{n-n'}}{b^{n-n'}}\right]$$
where the final factor on the right hand side is equal to an odd integer.  This requires that
$b^{n-n'}$ divide $2k'+1$ since $a$ and $b$ are relatively prime.  For each choice of $n' <n$,
every $2k'+1$ which is a multiple of $b^{n-n'}$ satisfies this property.  This gives an infinite
set for each $n'<n$, but by the previous $n=n'$ argument, only finitely many of them for each $n'$
can be included in the set $\Gamma$.  Since $n$ is a fixed finite integer, we still have a finite
collection $\Gamma$.

\item (\textit{$n'>n \geq 0$})\; We seek $\gamma'$ such that $$\gamma-\gamma' = \frac14
\left(\frac{b}{a} \right)^{n'} \left[ \frac{a^{n'-n}(2k+1)}{b^{n'-n}}-(2k'+1) \right],$$ where the
last expression is an odd integer.  Similar to the above argument, this requires that $b^{n'-n}$
divides $2k+1$.  Since $k$ and $n$ are fixed with $\gamma$, there can only be finitely many $n'>n$
such that $b^{n'-n}$ divides $2k+1$.  Then for each such $n'$, there can only be a finite selection
of values for $k'$ to include in $\Gamma$.

\end{enumerate}

This completes the proof that any set $\Gamma$ containing $0$ and elements of $\mathcal{O}$ such
that the difference between any two elements again is in $\mathcal{O}$ must be a finite set in the
case where $a$ is even and $b$ is odd.

\item \textbf{$a$ odd, $b$ even}\;  Now, let $\lambda = \frac{a}{b}$, where $b$ is even.  Since
$a$ and $b$ are relatively prime, $a$ must be odd.  We will show that the set $$\Gamma = \{0\} \cup
\left\{\frac{b^i}{4a}\;:\; i \in \N \right\}$$ has the property that the exponentials
$\{e_{\gamma}\,:\, \gamma \in \Lambda\}$ are pairwise orthogonal in $L^2(X_{\lambda},
\nu_{\lambda})$. We first show that the nonzero elements are all in $\mathcal{O}$.
$$\frac{b^i}{4a} = \frac14 \left(\frac{b^i}{a^i} \right) a^{i-1}$$ Since $a^{i-1}$ is an odd
integer for any choice of $i \geq 1$, each nonzero element of $\Gamma$ is in $\mathcal{O}$.

Next, we verify that differences of distinct nonzero elements of $\Gamma$ are in $\mathcal{O}$. Let
$\gamma = \frac{b^i}{4a}$ and $\gamma'=\frac{b^j}{4a}$, and assume without loss of generality that
$i < j$.

\begin{eqnarray*}  \gamma-\gamma' &=& \frac{b^i}{4a}-\frac{b^j}{4a}\\ &=& \frac{b^i(1-b^{j-i})}{4a}\\
&=& \frac14 \left(\frac{b}{a}\right)^i (1-b^{j-i})a^{i-1} \end{eqnarray*}

Since $1 \leq i < j$, the expression $(1-b^{j-i})a^{i-1}$ is an odd integer.  This proves that
$\gamma-\gamma' \in \mathcal{O}$.

\end{enumerate} \qed
\end{proof}

%%%%%%%%%%%%%%%%%%%%%%%%%%%%%%%%%%%%%%%%%%%%%%%%%%%%%%
\section{Experimental evidence for non-total sets}\label{sec:total}
%%%%%%%%%%%%%%%%%%%%%%%%%%%%%%%%%%%%%%%%%%%%%%%%%%%%%

Next, we consider whether any of our collections of orthogonal exponentials form orthonormal bases.
We can use Equation (\ref{eqn:product}) and Parseval's identity for orthonormal bases to determine
whether a collection of orthogonal exponentials is total.  From Parseval, we know that if
$\{e_{\gamma}\;:\; \gamma \in \Gamma\}$ is an ONB, then for any $f \in  L^2(X_{\frac34}, \nu_{\frac34})$,
we have  $$\|f\|_{\nu_{\frac34}}^2 = \sum_{\gamma \in \Gamma}|\langle f, e_{\gamma} \rangle|^2.$$
If we apply this to an exponential function, we find

\begin{eqnarray} \|e_t\|^2_{\nu_{\frac34}} &=& \sum_{\gamma \in \Gamma} \left|\langle e_t,
e_{\gamma}\rangle \right|^2 \nonumber\\ &=& \sum_{\gamma \in \Gamma}
|\widehat{\nu_{\frac34}}(t-\gamma)|^2 \qquad \mbox{by Equation (\ref{eqn:product}).}\label{onb}
\end{eqnarray}

Using Stone-Weierstrass to show the density of exponentials, we find that our collection is an ONB
if and only if the expression (\ref{onb}) is a function of $t$ identically equal to $1$.

\begin{equation}\label{eqn:total} \sum_{\gamma \in
\Gamma}[\widehat{\nu_{\frac34}}(t-\gamma)]^2 = \sum_{\gamma \in \Gamma}\prod_{k=0}^{\infty}
\cos^2\left(2\pi \Bigr(\frac34\Bigr)^k (t-\gamma)\right) \equiv 1 . \end{equation}

\begin{figure} \centering
\includegraphics{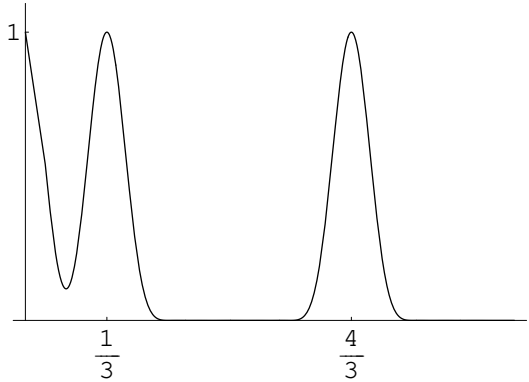}
\caption{$\Lambda_1$:  Equation (\ref{onb}) for $t\in [0,2]$. \label{lambda1}}
\end{figure}

\begin{figure} \centering
\includegraphics{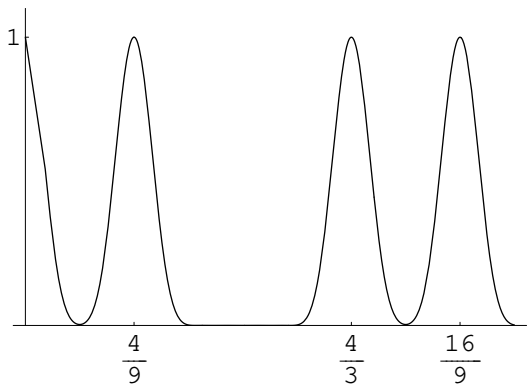}
\caption{$(\Lambda_1\cup \Lambda_2) \backslash \{1/3\}$:  Equation (\ref{onb}) for $t\in [0,2]$. We
gain a peak at $4/9$ but lose a peak at $1/3$.\label{lambda2}}
\end{figure}

\begin{figure} \centering
\includegraphics{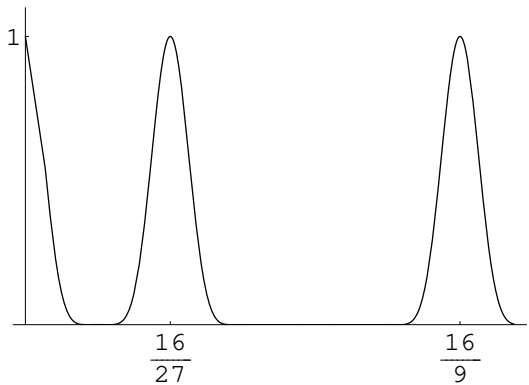}
\caption{($\Lambda_1\cup \Lambda_2\cup \Lambda_3) \backslash \{1/3, 4/3, 16/3\}$: Equation
(\ref{onb}) for $t\in [0,2]$.  We gain a peak at $16/27$ but lose at peaks at $4/9$ and $4/3$. \label{lambda3}}
\end{figure}

In Section \ref{sec:3/4}, we observed that we could remove one element from a $\Lambda_k$ set and
replace it with countably many other elements from another $\Lambda_j$ set.  For example, the set
$\Lambda_1 \cup \Lambda_2 \backslash \{1/3\}$ corresponds to a mutually orthogonal set of
exponentials.  However, when we graph the corresponding versions of Equation (\ref{onb}), we see
that the sum is far from being $1$.  See Figures \ref{lambda1}, \ref{lambda2}, and \ref{lambda3}
where we notice that as we omit a frequency, we lose a ``peak'' to $1$ at that frequency, but when
we include a frequency, we gain a ``peak'' to $1$.

\begin{figure} \centering
\includegraphics{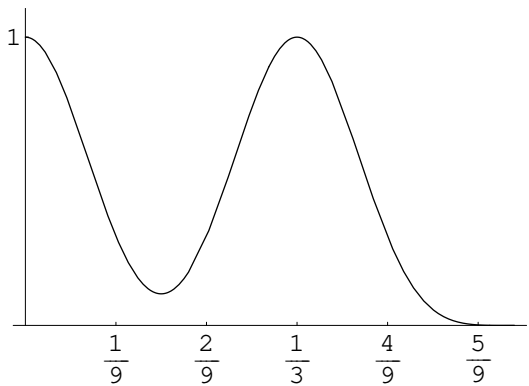}
\caption{The first $40$ terms of the sum
$\D \sum_{\ell\in \Gamma_1}
\prod_{n = 0}^{40}
\cos^2\Bigl(
2\pi \Bigl(\frac34\Bigr)^n
(\ell - t)
\Bigr).$
 The elements of $\Gamma_1$ are summed in
increasing order:  $0$, $1/3$, $4/3$, $5/3$, \ldots, $1045/3$.  \label{gamma1}}
\end{figure}

\begin{remark}An interesting corollary to Theorem \ref{thm:maximal} above is that when $\Gamma_1$
is used in the summation in Equation (\ref{eqn:total}), then as a function of $t$, this  $\Gamma_1$
summation is strictly positive on $R$.  (See Figure \ref{gamma1}.)
\end{remark}

We now turn to the sets $\Gamma_k$, which expand the $\Lambda_k$ sets in a different way from the
way we discussed above.  We have two reasons to suspect that none of the sets $\Lambda_k$ can be
used to construct ONBs.  First, numerical approximations from \textit{Mathematica} have provided
evidence that the sets $\Lambda_k$ are not total in $L^2(X_{\frac34}, \nu_{\frac34})$.  The graph in Figure
\ref{gamma1} shows that the expression (\ref{onb}) for $\Gamma_1$ is far from being identically $1$
after going out $40$ terms in the sum with $40$ terms in each product.  In contrast, the analogous
approximation of the sum and product for $\lambda = 1/4$ from \cite{JP1998}, for which there is an
orthonormal basis, appears to be identically $1$ after fewer than $40$ terms in each.

Second, Dutkay and Jorgensen have a conjecture \cite[Conjecture 6.1]{DJ2006b} which implies that in
the case $\lambda\in(1/2, 1)$, it is impossible to have an ONB for $L^2(X_{\lambda},
\nu_{\lambda})$.  In the conjecture, the existence of an ONB requires that a certain Hadamard
duality condition is fulfilled.   If the conjecture by Dutkay and Jorgensen is true, then none of
the sets $\Gamma_k$ (and therefore none of the sets $\Lambda_k$) for $\lambda\in(1/2, 1)$ can
possibly be ONBs.  Figure \ref{gamma1} provides graphical evidence for this conjecture in the case
of $\Gamma_1$ when $\lambda = \frac34$.

\section{Acknowledgements}
The authors would like to express their appreciation to Professor Kathy Merrill and Professor Judy
Packer for organizing the excellent conference \textit{Current Trends in Harmonic Analysis and Its
Applications: Wavelets and Frames} (affectionately known as C'estlarrybrate) in honor of Larry
Baggett.  Two of us attended this conference, with support from the conference funds, and our discussions there motivated some of the early stages of this work.  

The first named author was supported in part by a grant from the National Science Foundation.  Also, the first named author is pleased to acknowledge helpful discussions with Dorin Dutkay.  The second two named authors wish to express their appreciation of Grinnell College's contributing support for this work.

At the conclusion of our work on this paper, we were advised that existence results similar to
those in Theorems \ref{lambda_set_theorem} and \ref{thm:generallambda} were found simultaneously
and independently by Hu and Lau in \cite{HL2006}.

%\end{doublespace}

\end{document}